\documentclass[11pt, oneside]{article} 
\usepackage{latexsym}
\usepackage{amsfonts}
\usepackage{rawfonts}
\usepackage{graphicx}			
\usepackage{amsmath}										
\usepackage{amssymb}
\usepackage{geometry}                		
\input{prepictex}
\input{pictex}
\input{postpictex}

\newtheorem{defn}{Definition}
\newtheorem{thm}{Theorem}
\newtheorem{prop}[thm]{Proposition}
\newtheorem{cor}[thm]{Corollary}

\newtheorem{lem}[thm]{Lemma}

\newtheorem{pf}{Proof}
\newtheorem{rmk}{Remark}

\begin{document}

\title{On the uniqueness of conical K\"ahler-Einstein metrics}

\author{Long Li}

\maketitle{} 
\begin{abstract}
The purpose of this paper is to prove the uniqueness of conical K\"ahler-Einstein metrics, under the condition that the twisted $Ding$-functional is proper. This is a generalization of the author's previous work, and we shall first investigate the uniqueness of twisted K\"ahler-Einstein metrics, and then use these smooth perturbed solutions to approximate the actual conical one. Finally, it will bring the uniqueness of the limiting singular metric.
\end{abstract}
\section{Introduction}
The study of K\"ahler-Einstien equation plays an important role in nowadays geometric analysis after Yau's celebrated work on Calabi's Conjecture\cite{Yau}. The existence of these metrics were solved on the manifolds with non-positive first Chern class in Yau's paper in 1978\cite{Yau}. Then Matsushima\cite{Mat}, Futaki\cite{Fut} found some obstructions for the existence of K\"ahler-Einstein metrics on Fano manifolds($c_{1}(X) > 0$), and later Tian\cite{Tian83} gave the first sufficient condition for the existence in the Fano case. Many mathematicians contributed to this question in these 20 years, and most recently Chen, Donaldson and Sun's work\cite{CDS} proved that the existence of K\"ahler-Einstein metrics on Fano manifolds is equivalent to the K-stability condition. Their work based on an investigation on some conical K\"ahler metrics in a special H\"older space $\mathcal{C}^{2,\gamma,\beta}$ for 
$\gamma< \frac{1}{\beta} -1$, which was introduced by Donaldson\cite{Don12}. 
\\
\\
On the other hand, the uniqueness of K\"ahler-Einstein metrics on manifolds with $c_1(X) \leqslant 0$ was proved by Calabi in 1950s, and then Bando-Mabuchi\cite{BM} proved the uniqueness of K\"ahler-Einstein metrics on Fano manifolds($c_1(X)>0$) in 1987. It's natural to ask questions about uniqueness of conical or even more singular K\"ahler-Einstein metrics on Fano manifolds. For instance, BBGZ09(\cite{BBGZ}) proved the uniqueness of K\"ahler-Einstein metrics starting from a singular setting on a manifold $X$ without any non-trivial holomorphic vector fields. And Berman\cite{Berman12} proved the uniqueness of Donaldson's equation(a generalization of K\"ahler-Einstein equation in conical setting) based on a technical in Berndtsson's previous work\cite{Bo05} about curvature of $\log$ Begrman kernels. And then Berndtsson\cite{Bo11} proved some uniqueness results for K\"ahler-Einstein metrics with only $L^{\infty}$ potentials. 
\\
\\
The goal of this paper is to prove the uniqueness of $\mathcal{C}^{2,\gamma,\beta}$ conical K\"ahler-Einstein metrics in another way, based on a new technique developed in the author's previous work\cite{Li}. The main ingredients of this technique consist of Chen's $\mathcal{C}^{1,\bar{1}}$ geodesics and a generalization of Futaki's formula. In the first few sections, this technique will be extended to prove the uniqueness of the so-called twisted K\"ahler-Einstein metrics, i.e. a smooth K\"ahler metric $\omega$ satisfies
\[
Ric(\omega) = \omega + \theta
\]
where $\theta$ is some non-negative closed (1,1) form on $X$. In fact, assuming the correct cohomology condition(see section 2), we have 
\begin{thm} \label{utke1}
Suppose $\omega_0$ and $\omega_1$ are two solutions of twisted K\"ahler-Einstein equation with the same weight $\theta$, then there exists a holomorphic automorphism $F$ on $X$, such that $F^*(\omega_1) = \omega_0$ and $F^*(\theta) = \theta$, and this automorphism is induced from a holomorphic vector field $\mathcal{V}$. Moreover, if there is a point $p\in X$ such that the twister $\theta$ is strictly positive, then $\omega_0$ is actually fixed, i.e.
\[
\omega_1 = \omega_0
\]
on $X$.
\end{thm}
Next, before considering the singular metrics, we shall investigate the perturbed conical K\"ahler-Einstein equation first, i.e. we put  
\[
\theta = (1-\beta) dd^c \chi_{\epsilon}
\]
where the twister $\chi_{\epsilon} = \log (|s|^2 + \epsilon e^{\psi} )$ and $\psi$ is some smooth positively curved metric on the line bundle associated with the divisor $D$. Then we have the uniqueness of the solution of the following equations
\[
Ric(\omega_{\epsilon}) = \omega_{\epsilon} + dd^c \chi_{\epsilon}.
\]
If we can take the limit when $\epsilon\rightarrow 0$, then in principle, it will bring us the uniqueness of the conical K\"ahler-Einstein metrics $\omega_{\beta}$, i.e. $\omega_{\beta}$ satisfies 
\[
Ric(\omega_{\beta}) = \omega_{\beta} + (1-\beta)\delta_D
\]
where $\delta_D$ is the integration current of the divisor $D$. However, this is not true in general. In some situations, we can't not even find solutions of the perturbed K\"ahler-Einstein equations even if the conical K\"ahler-Einstein metric exists. Hence we need another condition to guarantee the existence of twisted K\"ahler-Einstein metrics, i.e. the properness of twisted $Ding$ functional. 

\begin{thm}\label{utke2}
Suppose the twisted $Ding$-functional $\mathcal{D}_{\beta}$ is proper, then there is only one $\mathcal{C}^{2,\alpha,\beta}$ solution $\omega_{\varphi_{\beta}}$ for the conical K\"ahler-Einstein equation with angle $\beta$ along the divisor $D$ on $X$.  
\end{thm}

One direct consequence of above theorem is the uniqueness of Donaldson's equation, i.e. 
\[
Ric(\omega_{\beta}) = \beta \omega_{\beta} + (1-\beta)\delta_{D}.
\]
for $0 < \beta < 1$. Here the K\"aher class is proportional to the anti-canonical class, and the twisted $Ding$ functional is automatically proper in this case\cite{Don12}. 
\\
\\
$\mathbf{Acknowledgement}$: I would like to express my gratitude to Prof. Xiuxiong Chen, who suggested me to do this problem and gave me consistently encouragement during the proof. I would also like to thank Prof. Song Sun, Prof. Dror Varolin for lots of helpful discussion. Finally, the suggestion from Chengjian Yao also helps me to make this paper more clear.

\section{Twisted K\"ahler-Einstein metrics}
Let $X$ be a compact complex K\"ahler manifold with $-K_X > 0$, and $S$ be a semi-positive $\mathbb{R}$-line bundle on $X$, i.e. there exists a smooth Hermitian metric $\psi$ on the line bundle $S$ such that 
\[
\theta = i\partial\bar{\partial}\psi \geqslant 0.
\]
Notice here $\theta$ is a globally defined closed $(1,1)$ form, and then we are going to consider the following twisted K\"ahler-Einstein equation
\begin{equation}
Ric(\omega) = \omega + \theta
\label{1}
\end{equation}
on $X$. Here we assume $-(K_X + S) > 0$, and a K\"ahler form $\omega$ can always be written as
\[
\omega = i\partial\bar{\partial}\phi
\]
where $\phi$ is a positively curved smooth metric on the $\mathbb{R}$-line bundle $-(K_X+S)$. And we will also use another notation when there is no confusion, i.e. 
\[
\omega_{\varphi} = \omega_0 + dd^c \varphi
\]
where $\omega_0$ is a fixed background K\"ahler metric in the same cohomology class, and $\varphi$ is the K\"ahler potential. We shall consider all such metric with $L^{\infty}$ potentials, i.e.
\[
|\phi' - \phi | < +\infty
\]
Let's denote these metrics as $PSH_{\infty}(-K_X-S)$, and we always assume $\int_X \omega^n =1$ in the following. Now suppose there exists two K\"ahler metrics $\phi_0$ and $\phi_1$ satisfying equation (\ref{1}), i.e.
\begin{equation}
\omega^n_{\phi_i} = \frac{e^{-\phi_i-\psi}}{\int_X e^{-\phi_i-\psi}}
\label{2}
\end{equation}
for $i =1,2$. We then claim they are the same up to a holomorphic automorphism, i.e. 
\begin{thm}\label{utke}
Suppose $\omega_1 = i\partial\bar{\partial}\phi_1$ and $\omega_2= i\partial\bar{\partial}\phi_2$ are two solutions of twisted K\"ahler-Einstein equation with the same weight $\theta$, then there exists a holomorphic automorphism $F$ on $X$, such that $F^*(\omega_2) = \omega_1$ and $F^*(\theta) = \theta$. Moreover, this automorphism is induced from a holomorphic vector field $\mathcal{V}$. 
\end{thm}

This is the twisted version of uniqueness theorem on K\"ahler-Einstein metrics on Fano manifolds. In order to investigate this equation from variational methods, we shaw introduce the twisted 
$Ding$-functional, whose critical point corresponds to the twisted K\"ahler-Einstein metric.
\begin{defn}
The twisted $Ding$-functional $\mathcal{D}$ is a functional defined on the space of all plurisubharmonic metrics $\phi\in PSH_{\infty}(-K_X-S)$, such that
\[
\mathcal{D}: = -\mathcal{E}+ \mathcal{F}_{\psi}
\]
where 
\[
\mathcal{E}(\phi): = \frac{1}{n+1}\int_X \Sigma_{j=0}^n (\phi - \phi_0) \omega_{\phi}^j \wedge \omega_{\phi_0}^{n-j}
\]
and 
\[
\mathcal{F}_{\psi}(\phi): = -\log \int_X e^{-\phi-\psi}.
\]
\end{defn}

\begin{rmk}
Since $\phi$ is a metric on the line bundle $ -( K_X+S )$ and $\psi$ is a metric on $S$, we see $\tau = \phi + \psi$ is a metric on $-K_X$, and $e^{-\tau}$ is a volume form on $X$.
\end{rmk}

Notice that the two functionals $\mathcal{E}$ and $\mathcal{F}$ will be changed in a new normalization, but $Ding$-functional is a normalization invariant. Now suppose there is a smooth curve $\phi_t$ in the space of metrics, and we can compute derivatives of twisted $Ding$-functional on it, i.e.
\[
\frac{\partial\mathcal{D}}{\partial t} = \int_X \phi' (\omega^n_{\phi} - \frac{e^{-\phi-\psi}}{\int_X e^{-\phi-\psi}}).
\]
Hence twisted K\"ahler-Einstein metric is its critical point, and
\[
\frac{\partial^2 \mathcal{D}}{\partial t^2} = \int_X (\phi'' - |\partial \phi'|^2_{g_{\phi}})\omega^n_{\phi} + (\int_X e^{-\tau})^{-1} \{\int_X (\phi'' - (\pi_{\perp}^{\tau} \phi')^2) e^{-\tau}   \}
\]
where $\tau = \phi +\psi$, and the orthogonal projection is $\pi^{\tau}_{\perp} u =  u - \int u e^{-\tau} / \int_X e^{-\tau}$. Now suppose $\phi_t, 0\leqslant t \leqslant 1$ is a smooth geodesic connecting $\phi_1$ and $\phi_2$, we see
\begin{equation}
 (\int_X e^{-\tau})\frac{\partial^2 \mathcal{D}}{\partial t^2} = \int_X (|\bar{\partial} \phi'|^2_{g_{\tau}} - (\pi_{\perp}^{\tau} \phi')^2) e^{-\tau} + \int_X (|\bar{\partial} \phi'|^2_{g_{\phi}} - 
 |\bar{\partial} \phi'|^2_{g_{\tau}})e^{-\tau} 
\label{3}
\end{equation}
where $g_{\tau}$ is the metric corresponding to the K\"ahler form $\omega_{\tau} = i\partial\bar{\partial}\tau$, and hence 
\[
\omega_{\tau} \geqslant \omega_{\phi}. 
\]
This implies the second term on the RHS of equation(\ref{3}) is non-negative, i.e. $ g^{\alpha\bar{\beta}}_{\phi}\phi'_{,\bar{\beta}}\phi'_{,\alpha} \geqslant  g^{\alpha\bar{\beta}}_{\tau}\phi'_{,\bar{\beta}}\phi'_{,\alpha} $, and the first term is non-negative from the {\em Futaki's formula} with respect to the weighted Laplacian operator \cite{Li}
\[
\Box_{\tau} = \bar{\partial}^*_{\tau}\bar{\partial}.
\]
Hence we proved the following
\begin{prop}
The twisted $Ding$-functional is convex along smooth geodesics. 
\end{prop}

We can further observe that on the smooth geodesic $\phi_t$, the twisted {\em Ding-functional} must keep to be a constant, i.e. $\partial^2 \mathcal{D} / \partial t^2 \equiv 0$ on the geodesic, and then the following two equations
\begin{equation}
\delta_{\tau}(t): = \int_X (|\bar{\partial} \phi'|^2_{g_{\tau}} - (\pi_{\perp}^{\tau} \phi')^2) e^{-\tau}  = 0
\label{5}
\end{equation}
\begin{equation}
k(t): = \int_X (|\bar{\partial} \phi'|^2_{g_{\phi}} -  |\bar{\partial} \phi'|^2_{g_{\tau}})e^{-\tau}  = 0
\label{6}
\end{equation}
hold simultaneously for each $0\leqslant t \leqslant 1$. The equation (\ref{5}) implies there exists a time independent holomorphic vector field $ V = X^{\alpha}\frac{\partial}{\partial z^{\alpha}}$, where
\[
X^{\alpha} = g_{\tau}^{\alpha\bar{\beta}} \phi'_{,\bar{\beta}}
\]
such that the Lie derivative
\[
\mathcal{L}_{\mathcal{V}} (\tau_t) = 0
\] 
where $\mathcal{V}_t = V - \frac{\partial}{\partial t}$, i.e. the induced automorphism $F$ preserves the metric as
 \[
 F^*(\omega_{\tau}(t)) = \omega_{\tau}(0) 
 \]
 for any $0\leqslant t\leqslant 1$. Then from the twisted K\"ahler-Einstein equation we see
 \[
 F^*(Ric(\omega_{\phi_1})) = F^*(\omega_{\phi_1} + \theta) = \omega_{\phi_0} + \theta = Ric(\omega_{\phi_0}),
  \]
hence $(F^*\omega_{\phi_1})^n = (\omega_{\phi_0})^n$ because they are in the same cohomology. By the uniqueness of Monge-Amp\`ere equation, we get
\[
F^*(\omega_{\phi_1}) = \omega_{\phi_0},
\]
the equation
\[
F^*(\theta) = \theta
\]
follows directly. Up to here, we proved theorem \ref{utke} under the assumption of smooth geodesics. Moreover, the equation (\ref{6}) implies the twisted K\"ahler-Einstein metric is really unique if $\theta$ is strictly positive. Next, we shall prove the theorem in the case when there is only $\mathcal{C}^{1,\bar{1}}$ geodesic connecting two twisted K\"ahler-Einstein metrics.
\\
\\
\begin{pf}[of theorem \ref{utke}]
Let's consider the $\epsilon$-approximation geodesics connecting $\phi_1$ and $\phi_2$, i.e. the solution of the following equation
\[
(\phi'' - |\partial \phi'|^2_{g_{\phi}})\det g_{\phi} = \epsilon \det h
\]
with boundary values 
\[
\phi(0,z) = \phi_1(z);\ \ \ \phi(1,z)=\phi_2(z).
\]
Now if we define $f = \phi'' - |\partial \phi'|^2_{g_{\phi}} > 0$, then we can see
\[
f\det g_{\phi} = \epsilon \det h
\]
and the second time derivative of twisted $Ding$-functional is 
\[
\frac{\partial^2 \mathcal{D}}{\partial t^2} = -\int_X f \det g_{\phi} + (\int_X e^{-\tau})^{-1} \{ \int_X ( f + \delta_{\tau} + k ) e^{-\tau}  \}.
\]
Notice that $\int_X f_{\tau} e^{-\tau}$, $\int_X \delta_{\tau} e^{-\tau}$ and $\int_X k e^{-\tau}$ are all non-negative, hence
\[
\int_{X\times I} ( f + \delta_{\tau} + k ) e^{-\tau} dt \leqslant \epsilon C
\]
for some uniform constant $C$ independent of $\epsilon$. If put $f_{\tau} = \phi'' - |\partial \phi'|^2_{g_{\tau}} = f + k >0 $, we see that for almost everywhere $t\in [0,1]$, there exists a constant 
$C(t)$ and a subsequence $\epsilon_j(t)$ such that 
\[
\int_X f_{\tau} e^{-\tau} (\epsilon_j) < C(t) \epsilon_j;\ \ \ \ \int_X \delta_{\tau} e^{-\tau}(\epsilon_j) < C(t) \epsilon_j
\]
from {Fatou's lemma}. Furthermore, all the uniform $\mathcal{C}^{1,\bar{1}}$ estimates hold for metrics $\tau(\epsilon)$ on the approximation geodesics, because 
\[
\tau_{\epsilon} = \phi_{\epsilon} + \psi
\]
and $\psi$ is a fixed smooth twister here. Hence by considering the weighted Laplacian operator $\Box_{\tau}$, the same argument in {\em \cite{Li} } implies there exists a time independent holomorphic vector field $V = X^{\alpha}\frac{\partial}{\partial z^{\alpha}}$, such that the induced automorphism $F$ of $X$ will preserve the metric, i.e. $F_t^*(\omega_{\tau}(t)) = \omega_{ \tau}(0)$, then as we argued before in smooth case, 
\[
F^*(\omega_{\phi_1}) = \omega_{\phi_0}
\]
and 
\[
F^*(\theta) = \theta
\]
for all $t\in[0,1]$.
\end{pf}

In prior, the holomorphic vector field $V = X^{\alpha}\frac{\partial}{\partial z^{\alpha}}$ obtained from above proof is computed with respect to the metric $\tau$, hence
\[
X^{\alpha} = g^{\alpha\bar{\beta}}_{\tau} \phi'_{,\bar{\beta}}.
\]
However, we know the condition 
\[
F^*_t (\omega_{\phi_t} )   = \omega_{\phi_0} + \theta - F^*_t \theta
\]
for each $t\in [0,1]$ from the proof. Hence the geodesic $g_{\phi_t}$ is smooth both in space and time directions, and the twisted $Ding$-functional $\mathcal{D}$ keeps being a constant along the geodesic $\phi_t$, i.e. 
\[
\frac{\partial\mathcal{D}}{\partial t} \equiv 0
\]
for all $t\in [0,1]$. Hence $\phi_t$ keeps to be the local minimizer of the twisted $\mathcal{D}$-functional along the whole curve, i.e. it satisfies 
\[
Ric(\omega_{\phi_t}) = \omega_{\phi_t} + \theta.
\] 
From the same argument as above we have 
\[
F^*_t(\omega_{\phi_t}) = \omega_{\phi(0)} ;\ \ \ F^*_t (\theta) = \theta
\]
for each $0\leqslant t \leqslant 1$. Then the curve generated by the one parameter group of automorphisms $F_t$ coincides with the geodesic $\omega_{\phi_t}$, and this is to say the 
$\mathcal{C}^{1,\bar{1}}$ geodesic is in fact smooth. Moreover, we can prove 
\begin{cor}
\label{pp}
Suppose there is one point $p\in X$, such that the closed (1,1) form $\theta$ is strictly positive, i.e. $\theta (p) > 0$, then the twisted K\"ahler-Einstein metric $\omega_1$ is actually unique. 
\end{cor}
\begin{pf}
Suppose we have two different twisted K\"ahler-Einstein metrics $\omega_1$ and $\omega_2$, then we can assume there is a smooth geodesic connecting them by the argument before the corollary. And this curve is the one parameter group of Automorphism $F_t$ generated by a nontrivial holomorphic vector field $V$. Now take the time derivative in the integral equation, we see
\[
F^*_t (\mathcal{L}_{\mathcal{V}}(\omega_{\phi_t}) ) =0 ,
\]
and hence 
\[
\partial (V \lrcorner \ \omega_{\phi_t}) = \omega'_{\phi_t} = \omega'_{\tau}
\]
which implies
\[
\partial (V \lrcorner \ \theta ) = 0
\]
Now notice the existence of twisted K\"ahler-Einstein metrics on the manifold $X$ implies the first betti number $b_1 = 0$, i.e. there is no nontrivial harmonic (0,1) form, hence $V \lrcorner \ \theta =0$, i.e.
\[
X^{\alpha}\theta_{\alpha\bar{\beta}} = 0
\]
and then we can write 
\[
X^{\alpha} = g_{\phi}^{\alpha\bar{\beta}}\phi'_{\bar{\beta}};
\]
for any time $0\leqslant t \leqslant 1$. Now in a neighborhood of the point $p\in U$, the (1,1) form $\theta$ keeps to be strictly positive, i.e.
\[
\theta > \epsilon \omega 
\] 
in $U$. But the equation $X^{\alpha}\theta_{\alpha\bar{\beta}} = 0$ implies the holomorphic vector field $V$ is identically zero in an open set $U$, hence it must be identically zero on $X$, which is a contradiction. 

\end{pf}

Finally by combining the results from theorem \ref{utke} and corollary \ref{pp}, the proof of theorem\ref{utke1} is finished.

\section{Smooth perturbation of conical K\"ahler-Einstein metrics}
We shall introduce conical K\"ahler-Einstein metrics first in this section, and consider the perturbation of these metrics. Let $D$ be a smooth divisor on the manifold $X$, such that the associated line bundle $S_D$ of this divisor is semi-positive. Suppose the $\mathbb{R}$-line bundle $-(K_X+(1-\beta)S_D)$ is strictly positive, where $0 <\beta<1$ is any real number. Then each $L^{\infty}$ strictly plurisubharmonic metric $\phi$ on this $\mathbb{R}$-line bundle associates with a K\"ahler form 
\[
\omega_{\phi} = i\partial\bar{\partial}\phi  > 0
\]   
as before.
\begin{defn}
A singular K\"ahler metric $\omega$ is a conical K\"ahler metric with angle $\beta$ along $D$ if the following conditions are satisfied:\\
\\
(i)\ $\omega$ is a closed positive (1,1) current on X, and is smooth on $X\backslash D$;\\
\\
(ii)\ for every point $p\in D$, there exists a constant $C$, such that in $D\cap U=\{z_1=0 \}$, where $U$ is a coordinate neighborhood of $p$, we have 
\[
C^{-1}\omega_{\beta} \leqslant \omega \leqslant C \omega_{\beta}
\]
where $\omega_{\beta}$ is a local model conical metric on $D\cap U$, i.e.
\[
\omega_{\beta}:= \sqrt{-1} (\frac{dz^1\wedge d\bar{z}^1}{|z^1|^{2-2\beta}} + \Sigma_{i=2}^n dz^i\wedge d\bar{z}^i)
\]
\end{defn}

From now on, we shall suppose $\omega_{\phi}$ is always a conical K\"ahler metric, and then we can talk about the K\"ahler-Einstein equation in this setting: we call $\omega_{\phi}$ is a {\em  conical K\"ahler-Einstein metric} if the following equation is satisfied in current sense
\begin{equation}
Ric(\omega_{\phi}) = \omega_{\phi} + (1-\beta)\delta_{D}
\label{7}
\end{equation}
where $\delta_D$ is the integration current of the divisor $D$, and it's equivalent to the following Monge-Amp\`ere equation
\begin{equation}
\omega_{\phi}^n = \frac{e^{-\phi} / |s|^{2-2\beta} }{\int_{X} e^{-\phi}/ |s|^{2-2\beta}},
\label{8}
\end{equation}
where $D=\{s=0 \}$. Notice that $e^{-\phi}/|s|^{2-2\beta}$ is a volume form on $X$ by the cohomology condition, since $\log |s|^2 $ corresponds to a plurisubharmonic metric on the line bundle $S$. 

Before going to this singular case, it's natural to consider the perturbed version of this equation, i.e. we shall consider the following smooth approximation equation
\begin{equation}
\omega_{\phi_{\epsilon}}^n = \frac{e^{-\phi_{\epsilon}}}{\mu_{\epsilon}(|s|^2 + \epsilon e^{\psi})^{1-\beta}}
\label{9}
\end{equation}
where $\psi$ is any semi-positively curved smooth metric on the line bundle $S$, and 
\[
\mu_{\epsilon} = \int_X \frac{e^{-\phi_{\epsilon}}}{(|s|^2 + \epsilon e^{\psi})^{1-\beta}}
\]
is the normalization constant. This is equivalent to the following geometric equality
\begin{equation}
\label{10}
Ric (\omega_{\phi_{\epsilon}}) = \omega_{\phi_{\epsilon}} + (1-\beta) \chi_{\epsilon}
\end{equation}
where 
\[
\chi_{\epsilon} =dd^c \log{( |s|^2 + \epsilon e^{\psi} )}
\]
Now we claim that $\chi_{\epsilon}$ is an non-negative closed (1,1) form as
\begin{lem}
For any $\epsilon > 0$, we have
\[
dd^c \log{( |s|^2 + \epsilon e^{\psi} )} \geqslant 0 .
\]
\end{lem}
\begin{pf}
\[
\partial\bar{\partial} \log {( |s|^2 + \epsilon e^{\psi} )} 
\]
\[
= \partial \{ \frac{sd\bar{s} + \epsilon e^{\psi} \bar{\partial}\psi}{|s|^2 + \epsilon e^{\psi}}  \}
\]
\[
= \frac{\epsilon e^{\psi}}{|s|^2 + \epsilon e^{\psi}} \{ \partial\bar{\partial}\psi + | ds - s \partial\psi |^2  \},
\]
and notice the above term is non-negative if the smooth metric $\psi$ is. 
\end{pf}

Now suppose we have a smooth solution $\omega_{\phi_{\epsilon}}$ of equation(\ref{9}), then theorem(\ref{utke1}) tells us
\begin{cor}
The smooth solution of the perturbed conical K\"ahler-Einstein equation $\omega_{\phi_{\epsilon}}$ is actually unique. 
\end{cor}
\begin{pf}
It's enough to prove there exist one point $p\in X$, such that $\chi_{\epsilon}(p) > 0$. First notice that this is the case if the metric $i\partial\bar{\partial}\psi$ is not completely degenerate. 
Otherwise, we have 
\[
\partial\bar{\partial}\psi \equiv 0
\]
on $X$. But meanwhile we should have  
\[
\bar{\partial}\psi = \frac{d\bar{s}}{\bar{s}}
\]
if $\chi_{\epsilon}$ vanishes. And then $i\partial\bar{\partial}\psi = \delta_D$, where $\delta_D$ is the integration current of $D$, which is a contradiction. 
\end{pf}

\section{Construction of the perturbed solution}

In this section, we shall discuss the existence of the solution of the perturbed K\"ahler-Einstein equation. In general, this is unknown even if the solution of equation (\ref{8}) exists. So,  here we need an extra assumption: the twisted $Ding$-functional is proper. \\
\\
In the following, we shall write the twisted $Ding$-functional as $\mathcal{D}_{\epsilon}$ with respect to the smooth twister $\chi_{\epsilon}$, and $\mathcal{D}_{\beta}$ with respect to the singular twister $(1-\beta)\delta_D$, i.e. we have
\[
\mathcal{D}_{\epsilon}(\phi) = -\mathcal{E}(\phi) - \log \int_X  \frac{e^{-\phi}}{(|s|^2 + \epsilon e^{\psi})^{1-\beta}}
\]
and 
\[
\mathcal{D}_{\beta}(\phi) = -\mathcal{E}(\phi) - \log \int_X  \frac{e^{-\phi}}{|s|^{2-2\beta}}.
\]
Notice that the critical point of the twisted $Ding$-functional $\mathcal{D}_{\epsilon}$ (or $\mathcal{D}_{\beta}$) is the solution of the twisted K\"ahler-Einstein equation with twister $\chi_{\epsilon}$ (or $(1-\beta)\delta_D$). Next, we shall introduce another important functional {\em Aubin's J} functional, i.e. 
\[
\mathcal{J}(\phi): =  \int_X \varphi \omega^n_0 - \mathcal{E}(\varphi)
\]
where $\varphi = \phi - \phi_0$. This functional is a kind of $W^{1,2}$ norm of the potential $\varphi$, and we can compare it with $Ding$-functional
\begin{defn}
The twisted $Ding$-functional $\mathcal{D}_{\epsilon}  (or\  \mathcal{D}_{\beta})$ is called proper if there exists some constant $a>0$ and $b$ such that
\[
\mathcal{D}_{.}(\phi)  \geqslant a \mathcal{J}(\phi) + b
\]
for all $\phi\in PSH(-K_X - S)$.
\end{defn}
\begin{rmk}
Notice that the properness of the twisted $Ding$-functional is in fact independent of the twister as long as the twister is smooth. This is because the twisted $Ding$-functionals are comparable for two different smooth twisters, i.e. let $\psi_1$ and $\psi_2$ be the corresponding metrics associated to the twisters $\theta_1$ and $\theta_2$ , then there exists a constant $C$ such that
\[
-C < \psi_1 - \psi_2 < C
\]
then the twisted functionals satisfy 
\[
\mathcal{F}_{\psi_2}(\phi) = - \log \int_X e^{- \phi - \psi_2} = - \log \int_X e^{-\phi - \psi_1 + (\psi_1 -\psi_2)}
\]
hence 
\[
|\mathcal{F}_{\psi_1} - \mathcal{F}_{\psi_2} | < C'
\]
for some uniform constant $C'$. And since $\mathcal{E}$ and $\mathcal{J}$ functionals are independent of the twister, the assertion follows. 
\end{rmk}
\begin{rmk}
For the special twister $\chi_{\epsilon}$, observe that the major term $\mathcal{F}$ in the twisted $Ding$-functional has the following relation 
\[
\mathcal{F}_{\epsilon} \leqslant \mathcal{F}_{\epsilon'}
\]
if $\epsilon \leqslant \epsilon'$, hence $\mathcal{D}_{\epsilon}$ is decreasing when $\epsilon$ becomes smaller. Moreover, if $\epsilon = 0$, it achieves its minimum , i.e. 
\[
\mathcal{D}_{\beta} \leqslant \mathcal{D}_{\epsilon}
\]
for all $\epsilon > 0$ small. Hence in the practice, we can simply require $\mathcal{D}_{\beta}$ to be proper.
\end{rmk}

Generally speaking, the purpose for introducing properness of $Ding$-functional is to solve the following continuity path 
\[
Ric(\omega_{\phi_t}) = t \omega_{\phi_t} + (1-t)\omega_0 + (1-\beta)\chi_{\epsilon}.
\]
where $\omega_0$ is a fixed smooth K\"ahler metric in the same cohomology with $\omega_{\phi}$, and $\chi_{\epsilon} = i\partial\bar{\partial}\psi_{\epsilon}$ for 
$\psi_{\epsilon} = \log (|s|^2 + \epsilon e^{\psi})$. This equation is solvable up to $t=1$ if the twisted $K-energy$ \cite{Berman12} is proper, and this condition is equivalent to the properness of the twisted $Ding$-functional from the argument of Berman\cite{Berman12}, so the existence of smooth perturbed solutions are guaranteed from here. However, this is not the end of the story since we want to make sure these perturbed solutions can approximate the conical one in some sense. 

In order to do this, we need to consider a kind of conical K\"ahler metrics with better regularity. Donaldson\cite{Don12} introduced a special H\"older space $\mathcal{C}^{2,\alpha,\beta}$ where $0<\alpha< \frac{1}{\beta} - 1$ for real valued functions on $X$, and we can define the so called $\mathcal{C}^{2,\alpha, \beta}$ conical K\"ahler metric\cite{Don12} by requiring that a local K\"ahler potential lies in $\mathcal{C}^{2,\alpha,\beta}$, i.e.  in a local coordinate chart near the divisor $D$, we can always write 
\[
\omega = i\partial\bar{\partial}(\varphi +\psi)
\]
where $\varphi \in \mathcal{C}^{2,\alpha,\beta}$, and $\psi$ is some smooth function. Be aware that this condition is stronger than the condition of conical K\"ahler metrics. Simply speaking, for fixed angle $\beta$, a $\mathcal{C}^{2,\alpha,\beta}$ conical K\"ahler metric is a conical K\"ahler metric with uniform $\mathcal{C}^{2,\alpha, \beta}$ norm.

Now let's take $\omega_{\varphi_{\beta}}$ to be a $\mathcal{C}^{2,\alpha,\beta}$ conical K\"ahler-Einstein metrics on $X$ and assume the twisted $Ding$-functional $\mathcal{D}_{\beta}$ is proper, then Chen, Donaldson and Sun's work\cite{CDS} provided a way to construct the following a family of sequence of perturbed solutions 
\[
\omega_{\phi_{\epsilon}}(t,.) = \omega_0 + i\partial\bar{\partial}\phi_{\epsilon},
\]
and it satisfies the following equation
\[
Ric(\omega_{\phi_{\epsilon}(t)}) = t\omega_{\phi_{\epsilon}(t)} + (1 - t)\omega_{\varphi_{\epsilon}} + (1-\beta)\chi_{\epsilon}.
\]
Notice that for $t=1$, $\omega_{\phi_{\epsilon}}(1, .) = \omega_{\phi_{\epsilon}}$ is exactly the solution of the twisted K\"ahler-Einstein equation with smooth twister $(1-\beta)\chi_{\epsilon}$, i.e.
\[
Ric(\omega_{\phi_\epsilon}) = \omega_{\phi_{\epsilon}} + (1-\beta)\chi_{\epsilon},
\]
and when $t=0$, the metric $\omega_{\phi_{\epsilon}}(0, .) = \psi_{\epsilon}$ will approximate the original metric $\omega_{\beta}$ because we require
\[
\omega_{\psi_{\epsilon}}^n = e^{-\varphi_{\epsilon} + h_{\omega_0}}\frac{1}{(|s|^2_h + \epsilon)^{1-\beta}}\omega_0^n
\]
where $\omega_{\varphi_{\epsilon}}$ is a small perturbation of the original metric such that $\omega_{\varphi_{\epsilon}} \rightarrow \omega_{\beta}$ in $C^{\gamma}(X)$. Moreover, we know that when $\epsilon\rightarrow 0$, the sequence of metrics $\phi_{\epsilon}(t, .)$ converges to a metric $\phi_0(t, .)$ globally in $C^{\gamma}$ and locally in $C^{3,\gamma}$ outside the divisor $D$, and the limiting metric $\phi_0 (t, .)$ will satisfy the following equation outside $D$
\begin{equation}
\label{approx}
\omega_{\phi_0}^n = e^{-t\phi_0 - (1-t)\varphi_{\beta} + h_{\omega_0}} \frac{1}{|s|_h^{2-2\beta}}\omega_0^n
\end{equation}
with 
\[
\phi_0(0, .) = \varphi_{\beta}.
\]
Observe that above equation is in fact equivalent to 
\begin{equation}
\label{approx2}
\omega_{\phi_0}^n =  e^{-t(\phi_0 - \varphi_{\beta})}\omega^n_{\varphi_{\beta}}
\end{equation}
with $t\in [0,1]$. Now we claim the curve of metrics $\phi_0(t, 0)$ is in fact fixed. The argument of the claim is similar with the end of the paper \cite{CDS}, and we shall recall it here for convenience of the reader
\begin{lem}
\label{fix}
The continuous family $\phi_0(t, .)$ is independent of the time $t$, and hence 
\[
\phi_0(1, .) = \varphi_{\beta}.
\]
\end{lem}
\begin{pf}
We shall argue like the end of the paper \cite{CDS}. First notice that the weighted Laplacian operator $\Delta_{\varphi_{\beta}}$ is continuous and invertible as a map 
\[
\Delta_{\varphi_{\beta}}: \mathcal{C}_0^{2,\gamma,\beta}\rightarrow \mathcal{C}^{,\gamma,\beta}
\]
for some $\gamma < \frac{1}{\beta} - 1$, and $\mathcal{C}_0^{2,\gamma,\beta}$ is the space of functions in $\mathcal{C}^{2,\gamma,\beta}$ with zero average. In fact, by Donaldson's H\"older estimate of conical metrics in \cite{Don12}, we can prove the first eigenvalue of this operator is strictly positive, i.e. 
\[
\Delta_{\varphi_{\beta}} > \lambda 
\]
for some constant $\lambda > 0$, and any $u\in \mathcal{C}_0^{2,\gamma,\beta}$. Now since the term $e^{-t(\phi_0 - \varphi_{\beta})}$ lies in $\mathcal{C}^{,\gamma, \beta}$, 
Implicit Function Theorem implies the existence of a continuous family of solution of the following equation 
\[
\omega_{\psi}^n = e^{-t(\phi_0 - \varphi_{\beta})}
\]
with $\psi(t, .)\in \mathcal{C}_0^{2,\gamma,\beta}$, for $t\in [0, \epsilon_0)$. Notice that $\psi(0, .) = \varphi_{\beta}$, hence $\psi(t, .)$ must coincide with $\phi_0(t, .)$ up to a constant in this short time. Then we can guarantee that $\phi_0(t, .) \in \mathcal{C}^{2,\gamma,\beta}$ for $t\in [0,\epsilon_0 )$.

Next we consider a constant continuity path $\varphi_{t}  \equiv \varphi_{\beta}$, which satisfies the equation 
\[
\omega_{\varphi_{t}}^n = e^{-t(\varphi_{t} - \varphi_{\beta})} \omega^n_{\varphi_{\beta}}.
\]
Now by Implicit Function Theorem again, we have a unique path of solution for $t < \min (\epsilon_0, \lambda/2)$, and hence 
\[
\phi_0(t, .) = \varphi_{\beta}, \ \ \forall\  0 \leqslant t < \min (\epsilon_0, \lambda/2).
\]
Finally we can repeat this procedure again and again, it will reach $t=1$, i.e.
\[
\phi_0(1, .) =\varphi_{\beta}.
\]
\end{pf}

Recall from our previous construction, we see the sequence of twisted K\"ahler-Einstein metric $\phi_{\epsilon}$ will converges to $\varphi_{\beta}$ in $C^{\gamma}$, hence we have
\begin{prop}
Suppose $\omega_{\varphi_{\beta}}$ is a $\mathcal{C}^{2,\alpha,\beta}$ solution for the conical K\"ahler-Einstein equation along the divisor $D$ on $X$, and the twisted $Ding$-functional is proper, we can construct a sequence of perturbed K\"ahler-Einstein metrics $\phi_{\epsilon}$, such that
\[
\phi_{\epsilon}\rightarrow \varphi_{\beta}
\]
globally in $C^{\gamma}$, and locally in $C^{k,\gamma}$ for some $k>3$ outside of $D$. 
\end{prop}

and this brings the proof our uniqueness.
\begin{pf}[of theorem \ref{utke2}]
Suppose $\omega_{\varphi_{\beta}}$ and $\omega'_{\varphi_{\beta}}$ are two different $\mathcal{C}^{2,\alpha,\beta}$ conical K\"ahler-Einstein metrics on $X$ with cone angle $\beta$ along the divisor $D$. From above argument, we can construct a sequence of smooth twisted K\"ahler-Einstein metrics $\phi_{\epsilon_j}$ to approximate $\varphi_{\beta}$ in some H\"older space 
$C^{\gamma}$. Then repeat the construction again, we can find a subsequence $\phi_{\epsilon_{n(j)}}$ to approximate $\varphi'_{\beta}$ by the actual uniqueness of $\phi_{\epsilon}$, and this implies $\phi_{\beta} = \phi'_{\beta}$ on $X$, which is a contradiction. 
\end{pf} 

In the end of the paper, a new point in the above proof should be mentioned. Unlike the end of the paper\cite{CDS}, the vanishing of tangential holomorphic vector fields is not necessary in the proof of lemma({\ref{fix}}), i.e. the properness of twisted $Ding$ functional has not been used here. So we can try to generalize this method by starting from the construction of some suitable approximation of conical K\"ahler-Einstein metrics without properness of $Ding$ functional.

\end{document}